\documentclass[12pt,a4paper]{article}

\usepackage{amssymb}
\usepackage{url}

\newtheorem{theorem}{Theorem}
\newtheorem{proposition}{Proposition}

\newtheorem{lemma}{Lemma}
\newtheorem{definition}{Definition}

\def\qed{\hfill$\square$}

\begin{document}

%\mainmatter  % start of an individual contribution

% first the title is needed
\title{Intertwining Laplace Transformations of Linear Partial
Differential Equations}

%\titlerunning{Intertwining Laplace Transformations of Linear PDEs}

\author{Elena~I.~Ganzha%
\thanks{This paper was written with partial financial support from the KSPU grant
``Forming scientific collective \emph{Physics of nano- and
microstructures}''.
 }
% }
%\institute{%
 \\ Krasnoyarsk State Pedagogical University,\\
ul.~Lebedevoi, 89, Krasnoyarsk,
660049, Russia \\
%\email{
{\tt eiganzha@mail.ru}
 }
%\date{May 14, 2013}

  \maketitle

\begin{abstract}
We propose a generalization of Laplace transformations to the case
of linear partial differential operators (LPDOs) of arbitrary order
in ${\mathbb R}^n$. Practically all previously proposed differential
transformations of LPDOs are particular cases of this transformation
(intertwining Laplace transformation, $\mathcal{ILT}$). We give a
complete algorithm of construction of $\mathcal{ILT}$ and describe
the classes of operators in ${\mathbb R}^n$ suitable for this
transformation.
%    MSC: ????.

\end{abstract}

\section{Introduction}\label{sec-intro}

In the past decade a number of publications
\cite{GLT,GanzhaDini,She12,She13b,TS,Tsarev2006,Tsarev2005,ts09,TT2009}
were devoted to application of various differential substitutions to
construction of algorithms for closed-form solution of linear
partial differential equations or systems of such equations. The
obvious drawback was just the vast diversity of such differential
substitutions, often considered as absolutely different in
properties and necessary tools for their study. As we show in this
paper practically all the aforementioned approaches can be naturally
unified into a very simple new class of \emph{Intertwining Laplace
Transformations}. We start with the classical Laplace cascade
method.
% of integration of a scalar second-order equation on the plane.
It is well known (see \cite{Darboux,gour-l,Kaptsov,GLT}) that
second-order linear hyperbolic equations on the plane, for example
\begin{equation}
 Lu = u_{xy} + a(x,y)u_{x} + b(x,y)u_y + c(x,y)u = 0 \, , \label{eq_1}
\end{equation}
admit the classical Laplace transformation  based on the following
equivalent forms of (\ref{eq_1}):
\begin{equation}
 \left[(D_x + b)(D_y + a) - h(x,y) \right] u = 0    \label{eq_2}
\end{equation}
or
\begin{equation}
  \left[(D_y + a)(D_x + b) - k(x,y) \right] u = 0 \, ,       \label{eq_3}
\end{equation}
where $D_x = \frac{\partial}{\partial x}$, $D_y =
\frac{\partial}{\partial y}$, $h = a_x + ab -c$, $k = b_y + ab -c$.
Equation (\ref{eq_2}) is equivalent to the first-order system
\begin{equation}
\left\{\begin{array}{l}
  (D_y + a)u =  u_1, \\
  (D_x + b)u_1 = hu.
\end{array}\right.     \label{eq_4}
\end{equation}
If $h \neq 0$, we can find $u$ from the second equation of the
system (\ref{eq_4}) and substituting it into the first equation of
the system (\ref{eq_4}) we obtain the transformed equation
%\begin{equation}
$
 L_1u_1 = (u_1)_{xy} + a_1(x,y)(u_1)_{x} + b(x,y)(u_1)_y +
    c_1(x,y)u_1 = 0
$.
%    \, .
%    \label{eq_5}
%\end{equation}
The operator $L_1$ is called the Laplace $X$-transformation of the
operator $L$. A simple calculation shows that $L$ and $L_1$ are
connected by an intertwining relation $(D_y +a_1)L = L_1(D_y+a)$.
 Analogously, if the invariant $k\neq 0$, we can define
the Laplace $Y$-transformation of the operator $L$ using
(\ref{eq_3}).

The transformations described above underlie the classical algorithm
for finding solutions of certain equations of the form (\ref{eq_1})
(the Laplace cascade method). Namely, applying for example the
Laplace $X$-transformation several times, in some cases, one can
obtain an equation of the form $
 \big(D_x + b\big)\big(D_y + \hat a\big)u =0 \, ,
$ which can be integrated in quadratures. Then with the help of the
inverse Laplace $Y$-transformation, its complete solution can be
used to obtain the complete closed form solution of the original
equation (\ref{eq_1}). See \cite{Darboux,gour-l} for more detail.

In \cite{GLT} we described a simple method (actually dating back to
Legendre, cf.~\cite{gour-l}) to apply the Laplace transformation to
the general second-order linear hyperbolic equations on the plane
\begin{equation}
Lu = u_{xx}  +  B(x,y)u_{xy}  + C(x,y)u_{yy}  +  D(x,y)u_x  +
E(x,y)u_y  +  F(x,y)u  = 0 \, . \label{eq_6}
\end{equation}
Any hyperbolic equation (\ref{eq_6}) can always be written in the
characteristic form
\begin{equation}
Lu = (X_1X_2 - H)u = 0 \, ,  \label{eq_7}
\end{equation}
where the coefficients of the operators $X_i = D_x +
\lambda_i(x,y)D_y + \alpha_i(x,y)$ and the function $H(x,y)$ are
constructively found from the coefficients of the original equation
(\ref{eq_6}) (see \cite{GLT,GanzhaDini}). Using (\ref{eq_7}) we see
that (\ref{eq_6}) is equivalent to the system
\begin{equation}
\left\{\begin{array}{l}
  X_2u =  v, \\
  X_1v = Hu .
\end{array}\right. \label{eq_8}
\end{equation}
If $H = h(x,y)$ is a nonzero function, then from the second equation
%of the system
(\ref{eq_8}) we find  $u = H^{-1}X_1 v $. Substituting
it into the first equation of
%the system
(\ref{eq_8}) we obtain the transformed equation
\begin{equation}
L_1v = (X_2X_1 + \omega X_1- H)v = 0 \, ,  \label{eq_9}
\end{equation}
where
\begin{equation}
 \omega = - [X_2, H]H^{-1} \, ,\label{eq_10}
\end{equation}
where $[\ ,\ ]$ denotes the usual operator commutation. The operator
$L_1$ is the result of the so-called $X_1$-Laplace transformation
applied to the operator $L$. One can easily check that the operators
$L$ and $L_1$ are connected by the \emph{intertwining relation}
\begin{equation}
 M_1  L = L_1  M \, ,\label{eq_11}
\end{equation}
where $M = X_2$, $M_1 = X_2 + \omega$.

In \cite{GanzhaDini} we have described a generalization of the
Laplace transformation to second-order linear partial differential
operators in ${\mathbb R}^3$ (and, generally, in ${\mathbb R}^n$)
with the principal symbol decomposable into the product of two
linear factors. This generalization is based on the fact, that such
operators can be always represented in the form (\ref{eq_7}), but
the coefficients $\alpha_i$ of the operators $X_i$ and the term $H$
are elements of the noncommutative ring of differential operators
$\textbf{F}[D_z]$.

In dealing with these operators it is reasonable to use the
algebraic construction of the noncommutative Ore field
(\cite{Ore1,Ore2}) of formal ratios of differential operators
$\textbf{F}(D_z) = \{R \,|\, R = P^{-1}Q; P,Q \in\textbf{F}[D_z]\}$
(where  $\textbf{F}$ is some differential field of functions) with
the equivalence relation: $P^{-1}Q \sim K^{-1}N$, if there exist
$S,T \in\textbf{F}[D_z]$, $S\neq0$ and $T\neq0$ such that $SP =TK$
and $SQ = TN$.

More generally, any noncommutative ring $\textbf{K}$ for which the
Ore conditions are satisfied (see below) is isomorphically embedded
into the field $\textbf{T} = \{R| R = P^{-1}Q; P,Q \in\textbf{K}\}$
with the equivalence relation given above. The Ore conditions on the
original noncommutative ring $\textbf{K}$ are as follows:
\begin{enumerate}
\item $\textbf{K}$ contains no zero divisors; i.e., if $AB =0$ for some
 $ A,B \in\textbf{K}$, then $A = 0$ or $B = 0$;
\item $\forall A,B \in\textbf{K}$, $A \neq0$, $B\neq0$,
$\exists P \neq0$ and $Q\neq0$ such that $PA = QB$, and $\exists M
\neq0$ and $N\neq0$ such that $AM=BN$.
\end{enumerate}
It can be easily seen that the ring $\textbf{F}[D_z]$ (and rings of
operators with partial derivatives such as
$\textbf{F}[D_x,D_y,D_z]$) meet the Ore conditions; therefore,
$\textbf{F}[D_z]$ is isomorphically embedded into the above-defined
skew Ore field $\textbf{F}(D_z)$. This allows us to apply to the
coefficients of the operator (\ref{eq_7}) all arithmetical
operations, taking into account the property of noncommutativity.
This will always result in a differential operator of the same type
with coefficients in $\textbf{F}(D_z)$. The skew field
$\textbf{F}(D_z)$ has external derivations $D_x$ and $D_y$, which
can evidently be extended from the initial ring $\textbf{F}[D_z]$.
In the field $\textbf{F}(D_z)$, the order of an element is
determined correctly by the formula $\mbox{ord\,}(P^{-1}Q) =
\mbox{ord\,}(Q) - \mbox{ord\,}(P)$.

In \cite{GanzhaDini} we showed that formulas
(\ref{eq_7})--(\ref{eq_11}) also hold for second-order operators in
${\mathbb R}^n$. But constructive results were obtained in
\cite{GanzhaDini} only for operators with decomposable principal
symbol. For second-order equations in ${\mathbb R}^3$ this
decomposability means that the principal symbol taken as
second-order polynomial in formal commutative variables $\xi_i =
\frac{\partial}{\partial x_i}$ is decomposed into the product of two
polynomials that are linear with respect to $\xi_i$. This
restriction means that the operator $H$ in (\ref{eq_7}) is a
first-order operator with respect to $D_z$ only.

In the present paper we give the definition of a natural
generalization of the classical Laplace transformation for arbitrary
operators in ${\mathbb R}^n$ without any restriction on
decomposability of the principal symbol. We will call below such
generalization an \emph{Intertwining Laplace Transformation}
($\mathcal{ILT}$). We prove some general properties of such
transformations, demonstrate its generality on a wide range of
examples and give the general algorithm of construction of
$\mathcal{ILT}$ in  ${\mathbb R}^n$.

The paper is organized as follows.
 We give the general definition of the $\mathcal{ILT}$ in Sect.~\ref{sec-def} and
the general algorithm for its construction in Sect.~\ref{sec-alg}.
 Section~\ref{sec-nonex} contains a result on non-existence of  $\mathcal{ILT}$
for a generic second-order operator in ${\mathbb R}^n$ for $n \geq
3$. The generality of the notion of $\mathcal{ILT}$ is demonstrated
in Sect.~\ref{sec-intertw} on many famous first-order differential
transformations of linear ordinary and partial differential
equations.
 In Sect.~\ref{sec-map} we discuss surjectivity and invertibility of $\mathcal{ILT}$.
Section~\ref{sec-conc} contains concluding remarks on possible
future developments, in particular the statement of a general result
on representability of arbitrary intertwining relation (\ref{eq_11})
with first-order intertwining operator $M$ and arbitrary linear
partial differential operator $L$ in ${\mathbb R}^n$ as
$\mathcal{ILT}$. The Appendix contains an important technical result
establishing a correspondence between existence of an intertwining
relation (\ref{eq_11}) and existence of the left least common
multiple of the operators $L$ and $M$ in the ring of linear partial
differential operators.

\section{Definition of Intertwining Laplace Transformations
($\mathcal{ILT}$)}\label{sec-def}

Let $L$ be a general linear differential operator of arbitrary order
in ${\mathbb R}^n$
% of order two:
% \vspace{-0.5em}
%\begin{equation}
% L = \sum_{i\leq j}^n  A_{ij}D_{x_i}D_{x_j} +
% \sum_{k=1}^n  A_{k}D_{x_k}
%  + A_0
%  \label{eq_12}
%\end{equation}
% \vspace{-0.5em}
%\noindent
with coefficients
%$A_\ast=A_\ast(x_1,\ldots,x_n)$
from some constructive differentially closed field of functions
$\textbf{F}$. Below to simplify notations we set $n=3$. However, all
results are true for arbitrary $n\geq 1$. Let $X_1$, $X_2$ be
arbitrary differential operators from the ring of linear partial
differential operators $\textbf{F}[D_x,D_y,D_z]$, then one can
always represent the operator $L$ in the following form
\begin{equation}
L = X_1X_2 - H \, ,  \label{eq_13}
\end{equation}
where $H = X_1X_2 - L$ is a differential operator in
$\textbf{F}[D_x,D_y,D_z]$ (in general of arbitrary order). We form
% the operator
\begin{equation}
L_1 = X_2X_1 + \omega X_1 - H \, ,  \label{eq_14}
\end{equation}
where
\begin{equation}
\omega = - [X_2, H]H^{-1}   \label{eq_omega}
\end{equation}
is a (pseudo)differential operator (an element of the skew Ore field
$\textbf{F}(D_x,D_y,D_z)$). It is easy to check that the
intertwining relation  (\ref{eq_11}) automatically holds with the
operators $M = X_2$ and $M_1 = X_2 + \omega$.

The formulas (\ref{eq_11}), (\ref{eq_13}), (\ref{eq_14}) hold in
$\textbf{F}(D_x,D_y,D_z)$. But it is difficult to use them for
transformation of solutions of the equation $Lu=0$ into solutions of
$L_1v=0$, the latter being a pseudodifferential equation in the
general case. So we introduce the definition of Intertwining Laplace
Transformation in which we impose the strong condition that $\omega$
should be a \emph{differential} operator (an element of the subring
$\textbf{F}[D_x,D_y,D_z]$):
\begin{definition} \label{def-ILT} We will say that the differential
operators $L$ and $L_1$ defined above by the formulas (\ref{eq_13})
and (\ref{eq_14}) with the condition $\omega = - [X_2, H]H^{-1}\in
\textbf{F}[D_{x_1},\ldots ,D_{x_n}]$, are connected by an
\emph{Intertwining Laplace Transformation} ($\mathcal{ILT}$).
 \end{definition}
\begin{lemma} If the operators  $L$ and $L_1$ are connected by
an $\mathcal{ILT}$ then their principal symbols coincide (even if
$\mbox{ord\,} H \geq \mbox{ord\,} L$).
 \end{lemma}
{\bf Proof.} %\begin{proof}
Since $\mbox{ord\,} \omega = \mbox{ord\,} X_2 -1$, the principal
symbols of the operators $M=X_2$ and $M_1=X_2 + \omega$ coincide.
Then (\ref{eq_11}) implies $\mbox{Sym\,} L = \mbox{Sym\,} L_1$. \qed
  % \end{proof}

%
It is easy to see that $\mathcal{ILT}$ is a generalization of the
classical Laplace transformation of  second-order operators in
${\mathbb R}^2$. However, it should be noted that even for dimension
two there exist other transformations different from the classical
Laplace transformation. They are defined with the help of the
intertwining relation (\ref{eq_11}) with some differential operators
$M$ and $M_1$. Such transformations were described in \cite{Darboux}
and will be considered in Sect.~\ref{subsec-lapdar}.

\section{Algorithm of Construction of $\mathcal{ILT}$ in ${\mathbb R}^n$}
\label{sec-alg}

As we will see below in Sections~\ref{sec-nonex},
\ref{subsec-gauge}-\ref{subsec-dini}, existence and construction of
intertwining relations (\ref{eq_11}) for a given operator $L$ is a
nontrivial problem. Even the ``functional dimension'' (number of
functions of maximal number of variables) of the set of all possible
pairs of operators $(L,M)$ admitting an intertwining relation
(\ref{eq_11}) is not known in the general case. In this Section we
give an algorithm which may be used to construct an arbitrary
$\mathcal{ILT}$ with first-order $M=X_2$. It should be noted however
that a given intertwining relations (\ref{eq_11}) may be represented
as an $\mathcal{ILT}$ in a non-unique way (see also
Sect.~\ref{sec-conc}).

\begin{lemma}\label{lemma1}
 For first-order intertwining operator $M = X_2$ the element
of the skew Ore field $\textbf{F}(D_x,D_y,D_z)$ $\omega = - [X_2,
H]H^{-1}$ is a differential operator (and consequently $L_1$ and
$M_1$ are differential operators) if and only if the operators $H$
and $X_2$ satisfy the relation
\begin{equation}
 HX_2 = (X_2 + \psi(x,y,z))H   \label{eq_15}
\end{equation}
with some function $\psi\in\textbf{F}$.
\end{lemma}
{\bf Proof.} %\begin{proof}
Let $\omega$ be a differential operator then we have $[H,X_2] =
\omega H$. Since $\mbox{ord\,} X_2 = 1$ we obtain $\mbox{ord\,}
\omega = 0$. Thus $\omega = \psi(x,y,z)$. This immediately implies
(\ref{eq_15}). The converse is obvious. \qed
  % \end{proof}

From Lemma~\ref{lemma1} immediately follows that if an
$\mathcal{ILT}$ connects operators $L$ and $L_1$ with a first-order
intertwining operator $X_2$ then $\omega$ is a function $\psi \in
\textbf{F}$. Below in such cases we will denote $\omega$ as $\psi$.

It is well known (the theorem on rectification of a vector field in
a neighborhood of each nonsingular point---a point where the vector
field is nonzero) that an arbitrary first-order operator $X_2$ may
be locally transformed to the form $X_2 = D_x + \alpha(x,y,z)$ with
$\alpha\in\textbf{F}$ by an appropriate (nonconstructive!)
coordinate transformation in a neighborhood of a generic point. For
this we need $\textbf{F}$ to be large enough to include the
necessary for this functions. In the new variables the relation
(\ref{eq_15}) has the form
\begin{equation}
 H(D_x + \alpha) = (D_x +\alpha + \psi)H \, .   \label{eq_16}
\end{equation}
Multiplying (\ref{eq_16}) on the left and on the right by some
functions $\mu(x,y,z)$ and $\rho(x,y,z)$ respectively we obtain
\begin{equation}\label{eq_16a}
 (\mu\,H\,\rho)(D_x + \alpha + \rho^{-1}\rho_x) =
 (D_x -\mu^{-1}\mu_x +\alpha + \psi)(\mu\, H\,\rho) \, .
\end{equation}
So the functions  $\mu$ and $\rho$ may be chosen in such a way that
(\ref{eq_16a}) will have the form of commutation relation
\begin{equation}
 \widetilde{H}D_x = D_x\widetilde{H} \label{eq_17}
\end{equation}
for the operator $\widetilde{H} = \mu\,H \,\rho$. Again we suppose
that $\textbf{F}$ is large enough to include $\rho$ and $\mu$. The
following Lemma may be easily proved by explicit computations.
\begin{lemma} The differential operator
$\widetilde{H}$ in ${\mathbb R}^n$ satisfies the commutation
relation (\ref{eq_17}) if and only if the coefficients of
$\widetilde{H}$ do not depend on $x$.
 \end{lemma}
Now we can formulate  \emph{the complete algorithm of construction
of arbitrary  $\mathcal{ILT}$ in ${\mathbb R}^n$}:
\begin{enumerate}
\item Take an operator $\widetilde{H}$ in
  ${\mathbb R}^n$ with coefficients not depending on a variable
  $x$.
\item Form the operator $H = \theta_1\widetilde{H}\theta_2$,
where $\theta_i$ are arbitrary functions  in  $\textbf{F}$. Using
the relation (\ref{eq_16a}) we find the functions $\alpha$ and
  $\psi$.
\item Make an arbitrary change of variables in ${\mathbb R}^n$ and
find the images of the operators $(D_x +\alpha)$, $H$ and the
function $\psi$  in the new variables. They are precisely the
operators $X_2$, $H$ and the function $\psi$ in (\ref{eq_15}).
% with $\omega = \psi$.
\item Taking  $L = X_1X_2 - H$ and $L_1 = X_2X_1 +
  \psi X_1 - H$ with an arbitrary operator
$X_1$ and $M = X_2$, $M_1 = X_2 + \psi$ we obtain a general example
of the  $\mathcal{ILT}$.
\end{enumerate}
\textbf{Remark.} Note that this algorithm is able to produce only
different examples of $\mathcal{ILT}$ with different operators $L$,
$L_1$, $M$. The problem of construction of $\mathcal{ILT}$ for a
given $L$ is very difficult in the general case and will not be
addressed here. Some particular methods of construction of
$\mathcal{ILT}$ for some classes of operators $L$ are given in
Section~\ref{sec-intertw}.

%Using the algorithm  we construct an example of second-order
%operators in ${\mathbb R}^3 = (x,y,z)$ connected by an
%$\mathcal{ILT}$. The operator  $X_1$ may be chosen arbitrary: we
{\bf Example} \emph{of second-order operators in ${\mathbb R}^3 =
\{(x,y,z)\}$ connected by an $\mathcal{ILT}$.}
 The operator  $X_1$
may be chosen arbitrarily: we take $X_1 = x^2D_y + xyD_z +1$.
Following the algorithm we take
 $H = xD_z^2x^2$ and find $X_2 = D_x + \frac{2}{x}$ and
$\omega = \psi = -\frac{3}{x}$. We omit the step 3, i.e. we will not
change the variables. Finally we obtain the operators
$$
L = X_1X_2 -H = x^2D_xD_y + xyD_zD_x -x^3D_z^2 + D_x + 2xD_y +2yD_z
+ {2}/{x} \, ,
$$
$$
L_1 = X_2X_1 + \psi X_1 -H = x^2D_xD_y + xyD_zD_x -x^3D_z^2 + D_x +
xD_y - {1}/{x} \, .
$$
The intertwining relation (\ref{eq_11}) has the form $ (D_x -
{1}/{x})  L = L_1 (D_x + {2}/{x}). $
% \end{example}

\section{On Non-existence of $\mathcal{ILT}$ for General
Second-Order Operators in ${\mathbb R}^n$}\label{sec-nonex}

\begin{theorem}\label{th-nonex}
 For a general second-order differential operator $L$
in ${\mathbb R}^n$ there are no $\mathcal{ILT}$ with first-order
operators $M$ for $n>2$.
\end{theorem}
{\bf Proof.} %\begin{proof}
Actually, following the algorithm we see that the number of
arbitrary functions of $n$ variables (we do not take into
consideration functions of smaller number of variables)
participating in the process of construction of $\mathcal{ILT}$ do
not exceed $2n + 3$ (two functions of $n$ variables $\theta_i$ on
step 2, $n$ functions on step 3 and $(n+1)$ coefficients of $X_1$).
On the other hand the number of the coefficients in a second-order
differential operator in ${\mathbb R}^n$ equals $\frac{n(n+1)}{2} +
n +1 > 2n + 3$, for $n>2$. Now the statement of the theorem is
obvious. \qed
  % \end{proof}

Note that the given estimate $2n + 3$ of the functional dimension of
the set of all $\mathcal{ILT}$ for second-order operators in
${\mathbb R}^n$ is just an upper bound, since different intermediate
data (operators $H$, $X_2$, functions $\theta_i$ etc. on the first
steps) may result in the same resulting operators $X_1$, $X_2$, $H$
in the final result of the algorithm. It would be interesting to
give a precise estimate of this functional dimension.

\section{Representation of Different Intertwining Relations as $\mathcal{ILT}$}
\label{sec-intertw}

In this section we show that many well known examples of
differential transformations of linear differential operators can be
represented as particular examples of the $\mathcal{ILT}$ introduced
in the previous Sections. We do this for:
\begin{enumerate}
\item gauge transformation $L\rightarrow\lambda^{-1}L\lambda$
where $\lambda \in \mathbf{F}$ is an arbitrary function;
\item differential substitutions for linear ordinary differential
operators and classical Darboux transformation for one-dimensional
Schr\"{o}dinger operator;
\item classical Laplace transformation and Darboux transformations for
$L = D_xD_y + a(x,y)D_x + b(x,y)D_y + c(x,y)$;
\item Euler-Darboux transformation (\cite{Kaptsov}) for  %high-order
operators in ${\mathbb R}^n$ of the form
$$
L = \sum_{i=0}^k  a_{i}(x)D_x^i + \sum_{|\alpha|\geq 0}^m
b_{\alpha}(y)D_y^{\alpha} \, ,
$$
where    $y = (y_1 , \ldots , y_{n-1})$, $\alpha = (\alpha_1, \ldots
, \alpha_{n-1})$, $ \displaystyle D_y^{\alpha} =
\frac{\partial^{|\alpha|}}{(\partial{y_1})^{\alpha _1}\cdots
(\partial{y_{n-1}})^{\alpha _{n-1}}}$;
\item  Darboux transformations for parabolic operators
$L = D_x^2 + a(x,y)D_x + b(x,y)D_y + c(x,y)$ on the plane;
\item Petr\'en transformation (\cite{Petren}) for higher-order operators
$$
L = \sum_{i=0}^{n-1}  A_i(x,y)D_xD_y^i +
 \sum_{i=0}^{n-1}  B_i(x,y)D_y^i \, ;
$$
\item Dini transformation (\cite{GanzhaDini,Tsarev2006}) for second-order operators in
${\mathbb R}^3$ with decomposable principal symbol.
\end{enumerate}
All aforementioned transformations are usually represented as
intertwining relations
\begin{equation}
 M_1  L = L_1  M \, . \label{int1}
\end{equation}
%with given operators $L$, $L_1$, $M$, $M_1$ and $\mbox{ord\,} L =
%\mbox{ord\,} L_1\geq 2$, $\mbox{ord\,} M = \mbox{ord\,} M_1\leq 1$.
From (\ref{int1}) we conclude that any solution $u$ of the equation
$Lu = 0$ is transformed into a solution $v = Mu$ of the transformed
equation $L_1v = 0$. Usually (\ref{int1}) is considered to be
fundamentally different from the classical Laplace transformation
described in Sect.~\ref{sec-intro}, since in many cases the mapping
$v = Mu$ of the solution space of the original equation $Lu = 0$ has
a nontrivial kernel (so is not invertible) unlike the classical
Laplace transformation (see Sect.~\ref{sec-map}).

We will use below the following precise definition of intertwining
relations:
\begin{definition}\label{def-intertw} Relation (\ref{int1}) with given
differential operators $L$, $L_1$, $M$, $M_1$ is called an
\emph{intertwining relation} between operators $L$ and $L_1$ with
\emph{intertwining operator} $M$ if the following conditions are
satisfied:
\begin{equation}
\mathrm{ord\,} L = \mathrm{ord\,} L_1 \, , \qquad
 \mathrm{ord\,} M = \mathrm{ord\,} M_1 \, , \label{ord-eq}
\end{equation}
\begin{equation}
\mathrm{Sym\,} L = \mathrm{Sym\,} L_1 \, . \label{Sym-eq}
\end{equation}
 \end{definition}
In Appendix we discuss the relation of the conditions
(\ref{ord-eq}), (\ref{Sym-eq}) to the existence of the left least
common multiple of the operators $L$ and $M$ in the ring of linear
partial differential operators.

First we establish a simple general result about representability of
(\ref{int1})
% with $\mbox{ord\,} M\leq 1$
as an $\mathcal{ILT}$.
\begin{proposition}\label{prop_X_1}
Let operators $L$, $L_1$, $M$, $M_1$ and its intertwining relation
$M_1L=L_1M$ be given. If there exists an operator $X_1$ which
satisfies the equation
\begin{equation}
 [X_1,X_2] - \omega X_1 = L -L_1 \, ,  \label{urX_1}
\end{equation}
with $X_2 = M$, $ \omega = M_1 - M$, then $L$ and $L_1$ are
connected by an  $\mathcal{ILT}$.
\end{proposition}
{\bf Proof.} %\begin{proof}
We should prove that all conditions of Def.~\ref{def-ILT} are
satisfied. Suppose $H=X_1X_2-L$ so $L=X_1X_2-H$ with $X_2=M$. We
should show that $[H,X_2]=\omega H$. Actually, $[H,X_2]=
[X_1X_2-L,X_2] = [X_1,X_2]X_2 - [L,X_2] = (\omega X_1 - L_1)X_2 +
X_2 L = \omega X_1X_2 - (X_2+\omega)L +X_2L = \omega H$. The
condition $L_1 = X_2X_1 + \omega X_1 - H$ follows automatically from
(\ref{urX_1}).
 \qed
  % \end{proof}

So (\ref{urX_1}) may be used to find $X_1$ in order to represent
(\ref{int1}) as an $\mathcal{ILT}$.
 In fact in many particular cases of intertwining
relations  we use another trick to find $X_1$ directly. This will be
explained in detail below.

\subsection{Gauge Transformation $L\rightarrow\lambda^{-1}L\lambda$,
\ \ $\lambda \in \mathbf{F}$} \label{subsec-gauge}

In this case the intertwining relation (\ref{int1}) is trivial
$\lambda^{-1}L = L_1\lambda^{-1}$, so $M = M_1 = X_2 =
\lambda^{-1}$, $\omega = 0$ and $X_1$ has to be found from the
condition (\ref{urX_1}), i.e.
$$
 X_1\lambda^{-1} - \lambda^{-1}X_1 = L  -\lambda^{-1}L\lambda \, .
$$
Obviously $X_1 = L\lambda + \varphi\lambda$ satisfies this equation
with arbitrary function $\varphi$. Then $H = X_1X_2 - L = \varphi$,
$L_1 = X_2X_1 - H = \lambda^{-1}L\lambda$, and all the required
relations of the $\mathcal{ILT}$ are satisfied.

\subsection{Differential Substitutions for Linear Ordinary Differential
Operators and Classical Darboux Transformation for One-Dimensional
Schr\"{o}dinger Operator}\label{subsec-lodo}

Here we consider first the so-called Loewy-Ore formal theory of
linear ordinary differential operators (LODO) which is described in
\cite{Tsarev1996}. For any two LODO $L$ and $M$ one can determine
their right greatest common divisor $rGCD(L,M)=G$, i.e.
$L=\widetilde L  G$, $M= \widetilde M  G$ (the order of $G$ is
maximal) and their left least common multiple $lLCM(L,M)=K$, i.e.
$K=\overline M L= \overline L M$ (the order of $K$ is minimal). This
can be done using the (noncommutative) Euclid algorithm in
$\mathbf{F}[D_x]$.  We say that {\it the operator $L$ is transformed
into $L_1$ by an operator} $M$, and write
$L\stackrel{M}{\longrightarrow}{L_1}$, if $rGCD(L,M)=1$ and $K=
lLCM(L,M) = L_1  M = M_1   L$. In this case any solution of $Ly=0$
is mapped by $M$ into a solution $z=My$ of $L_1z=0$. Using the
extended Euclid algorithm one may constructively find an operator
$N$ such that ${L_1}\stackrel{N}{\longrightarrow}L$, $N  M = 1 \,
(\mathrm{mod\,} L)$. Operators $L$, $L_1$ are also called {\it
similar} or {\it of the same kind} (in the given differential field
{\bf F} of their coefficients). So for similar operators the problem
of solution of the corresponding equations $Ly = 0$, $L_1y =0$ are
equivalent.

It is easy to represent the transformation
$L\stackrel{M}{\longrightarrow}L_1$ described above as an
$\mathcal{ILT}$. We consider the case $\mbox{ord\,} M = 1$ only.
Obviously, using the Euclidean division we obtain $L = QM + R$,
where R is a function and $R\neq 0$, since $rGCD(L,M) = 1$. Thus if
we take $X_1 = Q$, $X_2 = M$, $H = -R$, $\psi={[H,X_2]}/{H}$, so we
obtain the $\mathcal{ILT}$ with $L = X_1X_2 - H$, $L_1 = X_2X_1 +
\psi X_1 - H$, $M_1 = X_2 + \psi$ and the intertwining relation
\begin{equation}
 M_1L = L_1M \, , \label{int1a}
\end{equation}
where $\mbox{ord\,} M = \mbox{ord\,} M_1 = 1$, $\mbox{ord\,} L =
\mbox{ord\,} L_1$. Both sides of (\ref{int1a}) coincide with
$lLCM(L,M)$ since it is unique up to a factor $\alpha \in {\bf F}$.

Now let us consider the one-dimensional Schr\"{o}dinger operator
$$
 L = -\frac{d^2}{dx^2} + u(x) \, .
$$
Let $\omega$ satisfy the equation $L\omega = 0$.
The function $\omega$ determines a factorization of $L$:
\begin{equation} \label{onefactor}
L = A^\top A, \ \ \ A = -\frac{d}{dx} + v, \ \ \ A^\top =
\frac{d}{dx}+v, \ \ \ v = \frac{\omega_x}{\omega} \, .
\end{equation}
%Indeed we have
%$$
%A^\top A = \left(\frac{d}{dx} + v\right)\left(-\frac{d}{dx} +
%v\right) = -\frac{d^2}{dx^2} + v^2 + v_x
%and the equation
%$$
%v_x + v^2 = u
%$$
%is equivalent to $L\omega = 0$. If $v$ is real-valued we have
%$A^\ast = A^\top. $

{\it The Darboux transformation} is simply swapping of $A^\top$ and
$A$:
$$
L = A^\top A \to \widetilde{L} = AA^\top \, ,
$$
or in terms of the potential $u$: $u = v^2 + v_x  \to \widetilde{u}
= v^2 - v_x = u - 2(\log \omega)_{xx}$. So we have the intertwining
relation $AL = \widetilde{L}A$. In order to represent this
transformation as the $\mathcal{ILT}$ we take for example $X_1 =
A^\top$, $X_2 = A$, $H = 0$, $\psi = 0$, then $L = X_1X_2 - H =
A^\top A $, $L_1 = X_2X_1 + \psi X_1 - H = AA^\top  =
\widetilde{L}$. There is another possibility with $H \neq 0$.
Namely, we can take $X_1 = A^\top  + 1$, $X_2 = A$, $H = A$. Then
$[H,X_2] = [A,A] = 0$, thus $\psi = 0$, and we have $L_1 = X_2X_1 +
\psi X_1 - H = A(A^\top  + 1) - A = AA^\top = \widetilde{L}$.

\subsection{Classical Laplace Transformations and Darboux
Transformations for $L = D_xD_y + a(x,y)D_x + b(x,y)D_y +
c(x,y)$}\label{subsec-lapdar}

The  classical Laplace transformations for these operators are
obviously a particular case of the $\mathcal{ILT}$ (see
Sect.~\ref{sec-intro}). Another type of differential transformation
for this class of operators on the plane was studied by Darboux
\cite{Darboux}. Such Darboux transformation of order one is
constructed using a solution $u_1$ of the original equation $Lu =
0$. Darboux takes $M = D_x + \mu$, $\mu = -{(u_1)_x}/{u_1}$ or $M =
D_y + \nu$, $\nu = -{(u_1)_y}/{u_1}$, so that $Mu_1 = 0$, and proves
that there exists an operator $\widetilde L = D_xD_y + \widetilde
a(x,y)D_x + \widetilde b(x,y)D_y + \widetilde c(x,y)$ and $M_1 = D_x
+ \mu -{\alpha_x}/{\alpha}$, $\alpha = b + {(u_1)_x}/{u_1}$ (with
obvious modification for  $M = D_y -{(u_1)_y}/{u_1}$), which satisfy
$M_1L=\widetilde LM$. In order to represent this as an
$\mathcal{ILT}$ one should solve (\ref{urX_1}) for the unknown
operator $X_1$ with $\psi = -{\alpha_x}/{\alpha}$, $X_2 = D_x -
{(u_1)_x}/{u_1}$ and given $L$, $\widetilde L$. But we use another
trick taking into consideration the fact that the following system
\begin{equation}
\left\{\begin{array}{l}
  Lu = 0, \\
  X_2u = 0 .
\end{array}\right. \label{LX_2}
\end{equation}
has a nontrivial solution $u_1(x,y)$. We follow the usual way of
reducing (\ref{LX_2}) to involutive form, simplifying the first
equation of (\ref{LX_2}) using the second one:

$L = D_xD_y + a(x,y)D_x + b(x,y)D_y + c(x,y)$ $\longrightarrow$ $L -
(D_y + a(x,y))X_2 = \alpha(D_y  -\frac{(u_1)_y}{u_1})$, where
$\alpha = b + \frac{(u_1)_x}{u_1}$. We arrive at the system
\begin{equation}
\left\{\begin{array}{l}
  Hu = 0, \\
  X_2u = 0 ,
\end{array}\right. \label{HX_2}
\end{equation}
with $H = - \alpha(D_y  -\frac{(u_1)_y}{u_1})$. The system
(\ref{HX_2}) is obviously involutive and has one-dimensional
solution space generated by $u_1$ (see the standard techniques of
the Riquier-Janet theory for example in \cite{riq,jan,schwarz}).
This suggests to take $X_1 = D_y + a$ and $H$ given above. As one
can easily check the relation (\ref{urX_1}) with $L_1=\widetilde L$
is satisfied. Applying Prop.~\ref{prop_X_1} we come to the desired
representation of this Darboux transformation as an $\mathcal{ILT}$.

Darboux also studied transformations with higher-order operators
$M$, $M_1$ (as compositions of first order Darboux transformations).
We will limit ourselves to first order transformations and consider
the case $M = D_x + q(x,y)D_y + r(x,y)$ (with $q(x,y) \neq 0$). This
transformation is defined by two solutions $u_1$, $u_2 $ of the
original equation $Lu = 0$:
$$
%\begin{equation} \label{M1ord}
Mu = \left|
  \begin{array}{ccc}
    u   & u_y & u_x \\
    u_1 & (u_1)_y & (u_1)_x \\
    u_2 & (u_2)_y & (u_2)_x \\
  \end{array}
\right|\cdot\left|
  \begin{array}{cc}
    u_1 & (u_1)_y \\
    u_2 & (u_2)_y \\
  \end{array}
\right|^{-1},
%\end{equation}
$$
with the condition that $\left|
  \begin{array}{cc}
    u_1 & (u_1)_x \\
    u_2 & (u_2)_x \\
  \end{array}
\right| \neq 0$, $\left|
  \begin{array}{cc}
    u_1 & (u_1)_y \\
    u_2 & (u_2)_y \\
  \end{array}
\right| \neq 0$. We again form the system
\begin{equation}\label{(LMsys)}
\left\{\begin{array}{l}
  Lu = 0, \\
  Mu = 0
\end{array}\right.
\end{equation}
and reduce it to the involutive form
\begin{equation}
\left\{\begin{array}{l}
  Hu = 0, \\
  Mu = 0,
\end{array}\right. \label{(HMsys)}
\end{equation}
with $H = -(L -(D_y + a)M) = \gamma_2(x,y)D_y^2 + \gamma_1(x,y)D_y +
\gamma_0(x,y)$. This system is in involution since it has
two-dimensional solution space $\langle u_1, u_2\rangle$
(cf.~\cite{riq,jan,schwarz}). Obviously the commutator $[H,X_2]$
(where $X_2 = M$) is a second order operator $\theta_2 D_y^2 +
\theta_1D_y + \theta_0$ with some coefficients $\theta_i =
\theta_i(x,y)$ and  has the solution space $\langle c_1(x)u_1+c_2(x)
u_2\rangle$. The last operator has the same solution space as $H$,
so they are proportional: $[H,X_2] = \psi H$ for some function $\psi
= \psi(x,y)$. It gives us the necessary representation $L = X_1X_2
-H$ with $X_2 = M$, $X_1 = D_y +a$, $\psi =\theta_2/\gamma_2$ and
$H$ given above. Since the conditions of Prop.~\ref{prop-app} (see
Appendix) are obviously true for Darboux transformations we obtain
$\widetilde L = L_1= X_2X_1 + \psi X_1 -H$ and $M_1 = X_1 + \psi$.
Thus we again have represented this Darboux transformation with $M =
D_x + q(x,y)D_y + r(x,y)$ as an $\mathcal{ILT}$.

 Note that all considerations in this subsection are
valid for any hyperbolic operator $L$ with an arbitrary principal
symbol (\ref{eq_6}).

\subsection{Euler-Darboux Transformation  for  Higher-Order Operators
in ${\mathbb R}^n$ of the Form
 $
 L = \sum_{i=0}^k  a_{i}(x)D_x^i + \sum_{|\alpha|\geq 0}^m
 b_{\alpha}(y)D_y^{\alpha}
$ }\label{subsec-euler}
% where    $y = (y_1 , \ldots , y_{n-1})$, $\alpha = (\alpha_1,  \ldots ,
%\alpha_{n-1})$}

In this section we consider the linear partial differential equation
\begin{equation}\label{K1}
    Lu =Au + Bu =0 \, .
\end{equation}
Here $A$ is a differential operator w.r.t. the scalar variable $x$:
$
 A = \sum_{i=0}^k  a_{i}(x)D_x^i
$,
 and $B$ is a differential operator on the space of $n-1$
variables $y_1, \ldots y_{n-1}$: $
 B = \sum_{|\alpha|\geq 0}^m  b_{\alpha}(y)D_y^{\alpha}
 $,
where $y = (y_1 , \ldots , y_{n-1})$, $\alpha = (\alpha_1, \ldots,
\alpha_{n-1})$,
 $ \displaystyle D_y^{\alpha} =
\frac{\partial^{|\alpha|}}{(\partial{y_1})^{\alpha _1}\cdots
(\partial{y_{n-1}})^{\alpha _{n-1}}}$. We will denote  $E_{k,m}$ the
class of operators of the form (\ref{K1}). In \cite{Kaptsov} a
transformation of higher-order operators (\ref{K1}) was constructed
which generalizes the classical Euler (\cite{Euler}) and Darboux
(\cite{Darboux}) transformations  for  second order equations.
Following Kaptsov \cite{Kaptsov} we will call such transformation
Euler-Darboux transformation (${EDT}$). First we note that if
$h(x)$, $g(y)$ are solutions of the equations
\begin{equation}\label{K2}
    Ah =c \, h \, ,
\end{equation}
$$
Bg + c\,g =0, \ \ \ \ c\in \mathbb R \, ,
$$
then $u_1 = h(x)g(y)$ satisfies (\ref{K1}). The ${EDT}$ of the
operator $L$ is generated by its solution $u_1$. Namely (see
\cite{Kaptsov}) the differential substitution $w = hD_xh^{-1}u =
(D_x - {h_x}/{h})u$ maps solutions $u$ of (\ref{K1}) into solutions
$w$ of another equation $\widetilde{L}w = 0$ of the same class
$E_{k,m}$. This implies that the operators $L$ and $\widetilde{L}$
satisfy the intertwining relation $ M_1L = \widetilde{L}M$ with $M =
D_x - \frac{h_x}{h}$ and some first-order operator $M_1$. We again
do not solve (\ref{urX_1}) directly and use the same trick as in
Sect.~\ref{subsec-lapdar}. First we divide $A$ by $M$: $A = QM +
\varphi(x)$. Since $h(x)$ is a solution of (\ref{K2}) we see that
$\varphi (x) = const=c$ and $L = QM + c +B$. Now we can take $X_1 =
Q$, $X_2 = M$, $H = -(B +c)$ and obtain the necessary representation
for the operator $L$: $L = X_1X_2 - H$. We should check the
condition $[H,X_2] = \psi H$ for some function $\psi$. In fact
$[H,X_2] = [-B-c, D_x - \frac{h_x}{h}] = 0$ since the coefficients
of $B$ do not depend on $x$. Thus $\psi = 0$, and $L_1 = X_2X_1 - H
= MQ + B +c$. Using Prop.~\ref{prop-app} (see Appendix) we obtain
that $L_1=\widetilde{L}$.

\subsection{Darboux Transformations for Parabolic Operators
$L = D_x^2 + a(x,y)D_x + b(x,y)D_y + c(x,y)$}\label{subsec-parab}

We consider here the parabolic operator on the plane of the form
\begin{equation}\label{P1}
    L = D_x^2 + a(x,y)D_x + b(x,y)D_y + c(x,y), b(x,y)\neq 0 \, .
\end{equation}
In \cite{TS} the authors have proved that for any operator
(\ref{P1}) there exist infinitely many differential transformations
of the operator $L$ into the same form operators $\widetilde{L}$
which are defined by intertwining relation
\begin{equation}\label{P2}
    M_1L = \widetilde{L}M \, ,
\end{equation}
with operator $M$ of arbitrary order $k$ generated by some set of
independent solutions $z_1(x,y),\ldots , z_k(x,y) $ of the equation
$Lz = 0$. In contrast to the hyperbolic case considered in
Sect.~\ref{subsec-lapdar} there are no other differential
transformations similar to the classical  Laplace transformations.
We limit ourselves to the case of first-order operators $M = D_x +
q(x,y)D_y + r(x,y)$.

\underline{CASE A.} If $q\neq 0$ then the operator $M$ is defined by
conditions $Mz_1 = 0$, $Mz_2 = 0$ where $z_1$, $z_2$ are arbitrary
linearly independent solutions of $Lu = 0$ (\cite{TS}). Here we call
functions $z_1(x,y)$,  $z_2(x,y)$ \emph{linearly independent} if
they satisfy the following conditions:
$$
 \left|
  \begin{array}{cc}
    z_1 & (z_1)_x \\
    z_2 & (z_2)_x \\
  \end{array}
\right| \neq 0,  \qquad \left|
  \begin{array}{cc}
    z_1 & (z_1)_y \\
    z_2 & (z_2)_y \\
  \end{array}
\right| \neq 0 \, .
$$
 We
construct this operator $M$ as in Sect.~\ref{subsec-lapdar} using
the following analogue of the Wronskian formula:
$$
Mu = \left|
  \begin{array}{ccc}
    u   & u_y & u_x \\
    z_1 & (z_1)_y & (z_1)_x \\
    z_2 & (z_2)_y & (z_2)_x \\
  \end{array}
\right|\cdot\left|
  \begin{array}{cc}
    z_1 & (z_1)_y \\
    z_2 & (z_2)_y \\
  \end{array}
\right|^{-1} = ( D_x + q(x,y)D_y + r(x,y))u \, .
$$
In order to represent (\ref{P2}) as an $\mathcal{ILT}$ we take $X_2
= M = D_x + q(x,y)D_y + r(x,y)$. We  always can write the operator
(\ref{P1}) in the form $L = QX_2 + R$ where $Q = D_x -qD_y + (a-r)$
and $R = q^2D_y^2 + \alpha D_y + \beta $. Here $\alpha$, $\beta$ are
expressed in terms of the coefficients of the operators $L$ and $M$
and their derivatives. Setting $X_1 = Q$, $H = -R$ we come to the
required form $L = X_1X_2 - H$. The operators $H$ and $[H,X_2]$ are
second order operators containing only $D_y$ and satisfying the
condition $Hz_i = [H,X_2]z_i = 0, i =1,2$. The last condition
defines both operators up to a functional multiplier. So $[H,X_2] =
\psi (x,y)H$. This guarantees that the operators $L$ and $L_1 =
X_2X_1 + \psi X_1 - H$ are connected by $\mathcal{ILT}$. By
Prop.~\ref{prop-app} (see Appendix) we obtain that
$L_1=\widetilde{L}$.

\underline{CASE B.} If $q \equiv 0$ then the operator $M$   is
defined by one solution $z_1$ of $Lu = 0$ and should satisfy $Mz_1 =
0$ (\cite{TS}). The last condition implies that $M = D_x + r(x,y)$
where $r(x,y) = -{(z_1)_x}/{z_1}$. As before we set $X_2 = M = D_x +
r(x,y)$  and find the representation for the operator $L$: $L = QX_2
+ R$ with $Q = D_x + (a-r)$, $R = bD_y + (c - r_x + r^2 -ar)$.
Setting $X_1 = Q$, $H = -R$ we come to the required form $L = X_1X_2
- H$. The operators $H$ and $[H,X_2]$ are first order operators
containing only $D_y$ and satisfying the condition $Hz_1 =
[H,X_2]z_1 = 0$. This implies as in the Case~A that there exists a
function $\psi (x,y)$ such that $[H,X_2] =\psi(x,y)H$ and
$\widetilde{L} = L_1 = X_2X_1 + \psi X_1 - H$ with the help of
Prop.~\ref{prop-app}.

\subsection{Petr\'en Transformation (\cite{Petren}) for Higher-Order Operators
$L = \sum_{i=0}^{n-1}  A_i(x,y)D_xD_y^i + \sum_{i=0}^{n-1}
B_i(x,y)D_y^i$}\label{subsec-petren}

In \cite{LeRoux} a differential transformation for a class of
higher-order operators with two independent variables was proposed.
L.Petr\'en has extensively studied this transformation in her thesis
\cite{Petren}. Below we will call this transformation Petr\'en
transformation.

Petr\'en transformation applies to differential operators in
${\mathbb R}^2$ of the following form:
\begin{equation}\label{Pet1}
    L = \sum_{i=0}^{n-1}  A_i(x,y)D_xD_y^i + \sum_{i=0}^{n-1}  B_i(x,y)D_y^i \, .
\end{equation}
If we make a differential substitution $v = \alpha _0
D_y\alpha_0^{-1}u = (D_y - {(\alpha _0)_y}/{\alpha _0})u$ for any
solution $u$ of the equation $Lu = 0$ with the function $\alpha
_0(x,y) $ such that
\begin{equation}\label{Pet2}
    \sum_{i=0}^{n-1}  A_i(x,y)D_y^i\alpha _0 = 0 \, ,
\end{equation}
\begin{equation}\label{Pet3}
    L\alpha _0\neq 0 \, ,
\end{equation}
we obtain the transformed equation $\widetilde{L}v = 0$ with
operator $\widetilde{L}$ of the same type (\ref{Pet1}) (see
\cite{Petren}). As we have already seen this means that $L$ and
$\widetilde{L}$ are connected by an intertwining relation $
%\begin{equation}\label{Pet4}
    M_1L = \widetilde{L}M
%\end{equation}
$, where $M = D_y - {(\alpha _0)_y}/{\alpha _0}$ and $M_1$ is some
first-order differential operator. It will be shown in Appendix (see
Theorem 3 and the paragraph before it) that such operator exists and
its coefficients may be found constructively. We take $X_2 = M$. In
order to find $X_1$ we first write the operator $L$ in the form
\begin{equation}\label{Pet5}
  L = D_x\sum_{i=0}^{n-1}  A_i(x,y)D_y^i + \sum_{i=0}^{n-1}
   ( B_i(x,y) - (A_i(x,y))_x)D_y^i = D_x \widehat A + \widehat B \, .
\end{equation}
Using Euclidean division we can write
\begin{equation}\label{Pet6}
    \widehat A = \sum_{i=0}^{n-1}  A_i(x,y)D_y^i = QX_2 + q(x,y) \, ,
\end{equation}
where $Q$ is a differential operator of order $(n-2)$ and $q(x,y)$
is some function. Since $\widehat A\alpha _0 = 0$ by (\ref{Pet2})
and $X_2\alpha _0 = 0$ we have $q(x,y) \equiv 0$ in (\ref{Pet6}).
Analogously for the operator $\widehat B$ from (\ref{Pet5}) we have
\begin{equation}\label{Pet7}
   \widehat B = RX_2 + r(x,y) \, .
\end{equation}
Substituting (\ref{Pet6}) and (\ref{Pet7}) into (\ref{Pet5})  we
come to the form
\begin{equation}\label{Pet8}
    L = (D_xQ + R)X_2 + r(x,y) \, ,
\end{equation}
with $r(x,y) \neq 0$ by (\ref{Pet3}). Setting $X_1 = D_xQ + R$ and
$H = h(x,y) = -r(x,y)$ we obtain the required form $L = X_1X_2 -
h(x,y)$. The condition $[H,X_2] = \psi (x,y)H$ is obviously true
with $\psi(x,y) = -h_y/h$ since $h(x,y)$ is a function. By
Prop.~\ref{prop-app} from this follows that  $\widetilde{L} = L_1 =
X_2X_1 + \psi X_1 - H$. Thus Petr\'en transformation is represented
as an $\mathcal{ILT}$. We note that if $n = 2$ in (\ref{Pet1}) then
we will get the classical Laplace transformation with the
intertwining relation
$$
(D_y + A_0 - h_y/h)L_1 = L(D_y + A_0)
$$
where $h = (A_0)_x + A_0B_1 - B_0$, i.e. $h$ is the Laplace
invariant for (\ref{eq_1}).

\subsection{Dini transformation for second-order operators in
${\mathbb R}^3$ with decomposable principal
symbols}\label{subsec-dini}

In \cite{Tsarev2006}, an extension of the Dini transformation
(hereafter, simply ``Dini transformation'') was proposed that can be
applied to second-order operators in ${\mathbb R}^3$ with
decomposable principal symbols. Applied to such an operator
$L=X_1X_2 -H$ with first-order operators $X_1$, $X_2$, $H$, this
transformation takes a solution $u$ of the original equation $Lu =
0$ into the solution $v$ of the system
\begin{equation}\label{sisDini}
\left\{    \begin{array}{rl}
      (X_2 + \nu) u & =v, \\
      \widehat H u & = (X_1 + \mu)v.
    \end{array}\right.
\end{equation}
where $\widehat H = H + \mu X_2 + \nu X_1 + [X_1,\nu] + \mu\nu$ and
the functions $\mu, \nu \in \textbf{F}$ are to be chosen in such a
way that the old function $u$ can be eliminated from (\ref{sisDini})
obtaining a second-order equation $\widetilde{L}v=0$, rather than an
overdetermined system of equations for $v$ (which is the case for an
arbitrary system of the form (\ref{sisDini})). The latter means that
the commutator $\big[\widehat H,(X_2 + \nu)\big]$ can be expressed
in terms of the operators  $(X_2 + \nu)$, $\widehat H$ themselves:
\begin{equation}\label{uslDini}
  \big[\widehat H,(X_2 + \nu)\big]=\varkappa(x,y,z)\widehat H
  + \varrho(x,y,z) (X_2 + \nu) \, .
\end{equation}
So we obtain the equation $L_{Dini}v=0$ with the transformed
operator
\begin{equation}\label{LDini}
  L_{Dini}= (X_2 + \nu)(X_1 + \mu)
   -\widehat H +\varkappa(X_1 + \mu)  +\varrho.
\end{equation}
Condition (\ref{uslDini}) differs from our condition (\ref{eq_15})
by the presence of the second term in the right-hand side and
formally seems to be more general than (\ref{eq_15}). In
\cite{GanzhaDini} we have shown that in fact the conditions
(\ref{uslDini}) and (\ref{eq_15}) for existence of such $\mu$, $\nu$
are equivalent for the given operator $L$ and the resulting
transformed operators $L_{Dini}$ of the Dini transformation and
$L_1$ of $\mathcal{ILT}$ are also the same.
%This was achieved by a
%change of $\mu$ in (\ref{uslDini}). Details of the proof can be
%found in \cite{GanzhaDini}.
If we introduce $\widehat X_1 = X_1 + \mu$,  $\widehat X_2 = X_2 +
\nu$, then $L=\widehat X_1\widehat X_2 - \widehat H$;
(\ref{uslDini}) converts into
\begin{equation}\label{f42}
\big[\widehat H,\widehat X_2\big]=\varkappa(x,y,z)\widehat H +
\varrho(x,y,z) \widehat X_2,
\end{equation}
(\ref{LDini}) converts into
\begin{equation}\label{f43}
 L_{Dini}= \widehat X_2 \widehat X_1
-\widehat H +\varkappa\widehat X_1 +\varrho.
\end{equation}
Below we will always omit the hat sign over $\widehat X_i$ and
$\widehat H$ and write simply $X_i$, $H$, so for example $L_{Dini}=
X_2  X_1 - H +\varkappa X_1 +\varrho$, the same for (\ref{f42}) and
(\ref{f43}). Introducing an extra function $\alpha$ one  can write
the operator $L$ in the form
\begin{equation}\label{f44}
 L= X_1 X_2-H = (X_1 + \alpha)X_2 - \widetilde{H}
\end{equation}
with arbitrary $\alpha \in \mathbf{F}$ and $\widetilde{H} = H +
\alpha X_2$. By Def.~\ref{def-ILT} $L$ admits $\mathcal{ILT}$ if
there exists $\alpha$ such that the following condition is
satisfied:
\begin{equation}\label{f45}
 \big[\widetilde{H},X_2\big]=\psi\widetilde H
\end{equation}
with $\psi\in \mathbf{F}$.
\begin{proposition}\label{prop-Dini}
If some first-order operators $H$ and $X_2$ satisfy the condition
(\ref{f42}) with some functions $\varkappa$ and $\varrho$ then there
exists a function $\alpha$ such that $X_2$ and $\widetilde{H} = H +
\alpha X_2$ satisfy (\ref{f45}) with  $\psi=\varkappa$. The function
$\alpha$ is a solution of the equation
\begin{equation}\label{f-alfa}
\big[X_2, \alpha\big]+\varkappa \alpha = \varrho.
\end{equation}
\end{proposition}
{\bf Proof.} %\begin{proof}
Let $\alpha$ be an arbitrary function then obviously $ \big[H
+\alpha X_2,  X_2\big]= \big[H ,  X_2\big] + \big[\alpha,
X_2\big]X_2$. Using (\ref{f42}) for $\big[H ,  X_2\big]$ we see that
$ \big[H +\alpha X_2,  X_2\big]= \varkappa (H+\alpha X_2) +
\big(\varrho - \varkappa \alpha +\big[\alpha, X_2\big]\big)X_2 $.
Thus if $\alpha$ satisfies $\varrho - \varkappa \alpha +\big[\alpha,
X_2\big] = 0$ then (\ref{f45}) is satisfied with  $\psi=\varkappa$.
 \hfill \qed
\begin{theorem}\label{th-Dini}
Let $L$ be a second-order operator in ${\mathbb R}^n$ with
decomposable principal symbol and there exists its representation
$L=X_1X_2 -H$ with first-order operators  $X_i$, $H$ such that the
condition (\ref{f42}) is satisfied, i.e. $L$ admits Dini
transformation with the resulting operator $L_{Dini}$. Then there
exists a function $\alpha$ such that the operator $L$ represented in
the form $L= (X_1+\alpha)X_2 - (H +\alpha X_2) = \widetilde{X}_1X_2-
\widetilde{H}$ admits $\mathcal{ILT}$ with the resulting operator
$L_1 = X_2\widetilde{X}_1 + \psi \widetilde{X}_1 -\widetilde{H}
=L_{Dini}$.
\end{theorem}
{\bf Proof.} %\begin{proof}
From Prop.~\ref{prop-Dini} we obtain $ \big[\widetilde{H}, X_2\big]=
\big[H +\alpha X_2, X_2\big]= \varkappa (H+\alpha X_2)$ with
$\alpha$ satisfying (\ref{f-alfa}). So by Def.~\ref{def-ILT} we can
apply to $L$ the $\mathcal{ILT}$ and come to the transformed
operator $L_1 = X_2({X}_1 +\alpha) + \varkappa ({X}_1 +\alpha)
-\widetilde{H} = X_2X_1 + \varkappa X_1 - (H +\alpha X_2 - X_2
\alpha -\varkappa \alpha) =L_{Dini}$ since $\alpha$ satisfies
(\ref{f-alfa}).
 \hfill \qed

It should be noted that there is even a theorem in \cite{Tsarev2006}
stating that Dini transformation can be applied to any second-order
operators in ${\mathbb R}^3$  with a decomposable principal symbol
(i.e., appropriate $\mu$ and $\nu$ in operators $\widehat{X}_i$,
$\widehat{H}$ can always be found). Unfortunately, there is a grave
mistake in the proof of that theorem. In fact, Dini transformation
can be applied just to those operators to which the Intertwining
Laplace Transformation introduced here is applicable. As we have
shown in \cite{GanzhaDini}, this is not possible for arbitrary
second-order operators in ${\mathbb R}^3$ (see also
Theorem~\ref{th-nonex}).

\section{Mapping of the Solution Spaces for $\mathcal{ILT}$ and their Inverses}
\label{sec-map}

If we take the solution space  $S(L) = \{u| Lu = 0\}$ and any
intertwining relation (\ref{int1}) we obtain a linear mapping $M:
S(L)\longrightarrow S(L_1)$. Since in many cases considered in
Sect.~\ref{sec-intertw} and its subsections the operator $M$ has
solutions $z_i \in S(L)$, this mapping of the solution space has a
nontrivial kernel for such cases. As we have seen
\begin{equation}\label{LHX_2}
    \left\{\begin{array}{l}
  Lu = 0, \\
  Mu = 0,
\end{array}\right. \Longleftrightarrow \left\{\begin{array}{l}
  Hu = 0, \\
  Mu = 0,
\end{array}\right.
\end{equation}
here $X_2 = M$ and $H = X_1X_2 - L$. If $H$ and $X_2$ are
first-order operators with different principal symbols we conclude
from the condition $[H,X_2] = \psi (x,y)H$ that (\ref{LHX_2}) is
compatible and has one-dimensional solution space in ${\mathbb R}^2$
and infinite-dimensional solution space  in ${\mathbb R}^3$ (see the
basics of Riquier-Janet theory in \cite{riq,jan,schwarz}). In some
cases considered in Sect.~\ref{subsec-lapdar},~\ref{subsec-parab}
the operator $H$ had order two and the solution space of
(\ref{LHX_2}) was two-dimensional.

On the other hand it is easy to prove that in many cases the mapping
$M: S(L)\longrightarrow S(L_1)$ is surjective. This is true for all
cases studied in
Sect.~\ref{subsec-gauge},~\ref{subsec-lapdar}--\ref{subsec-dini} and
obviously not true for the one-dimensional Darboux transformation
(Sect.~\ref{subsec-lodo}). In fact we should prove in the
aforementioned nontrivial cases that for any $v\in S(L_1)$ there
exists $u \in S(L)$ such that $Mu =X_2u = v$. This means that the
following system
\begin{equation}\label{LX_2uv}
    \left\{\begin{array}{l}
  Lu = 0, \\
  X_2u = v
\end{array}\right.
\end{equation}
should have a solution iff $L_1 v = 0$. This system is obviously
equivalent to
$$\left\{\begin{array}{l}
  Hu = X_1v, \\
  X_2u = v.
\end{array}\right.
$$
In order to understand its compatibility conditions we again use the
Riquier-Janet theory \cite{riq,jan,schwarz}. Omitting the technical
details one gets the following result: existence of a solution of
this system is equivalent to the condition $X_2Hu - HX_2u = X_2X_1v
-  Hv$ (in fact the result of cross-differentiation of the equations
of this system if $ GCD(\mathrm{Sym\,} H , \mathrm{Sym\,} X_2)=1$).
Since $[H,X_2]u = \psi Hu = \psi X_1v$ we come to the equation
 $
X_2X_1v + \psi X_1v - Hv = L_1v = 0 % \, .
 $. Thus the system
(\ref{LX_2uv}) is compatible so the mapping $M: S(L)\longrightarrow
S(L_1)$ is surjective.

Let us write the transformed operator $L_1$ in the form
\begin{equation}\label{L_1X_2}
 L_1 = X_2X_1 + \psi X_1 - H = (X_2 + \psi )X_1 - H = \widetilde{X_2}X_1 - H \, ,
\end{equation}
where $\psi = [H,X_2]H^{-1}$. Then we can again apply to $L_1$ a
formal transformation  in the skew Ore field
$\textbf{F}(D_{x_1},\ldots ,D_{x_n})$ defined by the following
formulas:
\begin{equation}\label{L_1X_1sigma}
   (X_1 + \sigma)L_1 = \widetilde{L_1}X_1 \, ,
\end{equation}
$\sigma = - [X_1,H]H^{-1}$, $\widetilde{L_1} = X_1\widetilde{X_2} +
\sigma \widetilde{X_2} - H$. Note that $\sigma  \in
\textbf{F}(D_{x_1},\ldots ,D_{x_n})$ need not to be a differential
operator. It is easy to check that $\widetilde{L_1} = HLH^{-1}$.
Substituting it into (\ref{L_1X_1sigma}) we obtain $(X_1 +
\sigma)L_1 = HLH^{-1}X_1$ or $H^{-1}(X_1 + \sigma)L_1 = LH^{-1}X_1$.
% with $\sigma = - [X_1,H]H^{-1} = HX_1H^{-1} - X_1$.
Denoting $N = H^{-1}(X_1 + \sigma
)$, $N_1 = H^{-1}X_1$ we get the intertwining relation
 $
 NL_1 = LN_1
 $
which defines a formal transformation with $N, N_1 \in
\textbf{F}(D_{x_1},\ldots ,D_{x_n})$. This may be considered as a
pseudodifferential inverse of the $\mathcal{ILT}$ of the operator
$L$ into $L_1$. Note that we had to change  $X_2$ to
$\widetilde{X_2}$ in the representation of $L_1$ in order to obtain
this formal inverse. If one uses the representation $L_1 = X_2X_1 -
(H - \psi X_1) = X_2X_1 - \widetilde{H}$ and formally follows the
intertwining Laplace algorithm then the resulting operator will not
coincide with $L$.

Nevertheless one should note  that for some particular operators $L$
and mapping operators $M$ there exist differential operators $N$
mapping the solution space $S(L_1)$ onto $S(L)$ even if $M:
S(L)\longrightarrow S(L_1)$ has a nontrivial kernel. Examples of
such operators $L$, $L_1$, $M$, $N$ can be found in \cite{TS} where
existence of such differential operators $N$ was related to famous
nonlinear integrable equations for the coefficients of $L$.

\section{Conclusion}\label{sec-conc}

As we have demonstrated in the previous sections, the notion of
$\mathcal{ILT}$ unifies many differential transformations of linear
partial (and ordinary) differential equations previously considered
as fundamentally different.

The methods used in Sect.~\ref{sec-intertw} for representation of
various intertwining relations as $\mathcal{ILT}$ may be used to
prove the following general result:
\begin{proposition}\label{prop-conc}
Let $L$ be a linear differential operator of arbitrary order in
${\mathbb R}^n$ and the following intertwining relation
\begin{equation}\label{conc-int1}
  M_1L=L_1M
\end{equation}
holds for first-order operators $M$, $M_1$, and $\mbox{Sym\,} L =
\mbox{Sym\,} L_1$. Then $L$ is transformed by an $\mathcal{ILT}$ to
the operator $\alpha^{-1}L_1\alpha$ with $X_2=\alpha^{-1}M$, where
$\alpha$ is the coefficient at $D_{x_i}$ in $M$ (for any chosen
$i$).
\end{proposition}
The details of the proof and other developments on construction of
$\mathcal{ILT}$ for a given operator $L$ will be given elsewhere.
This proposition shows that the results of
Sections~\ref{sec-alg},~\ref{sec-nonex} may be actually formulated
for arbitrary intertwining relations with first-order operators $M$,
$M_1$.

A step into the direction of investigation of intertwining relations
with \emph{higher-order operators} $M$, $M_1$ may be potentially
obtained following the recent result \cite{She13b}, where the author
had proved that for the particular case of classical Laplace
operators (\ref{eq_1})  higher-order intertwining relations may be
represented as compositions of first-order $\mathcal{ILT}$.

Another important domain of applications for differential
transformations is the category of \emph{systems} of linear partial
differential equations, cf. for example
\cite{GLT,Tsarev2005,Athorne}. In fact, already
Le~Roux~\cite{LeRoux} had noted that it is much more natural to
study such transformations, since any differential substitution
$v=Mu$ transforms the solution space of a scalar equation $Lu=0$ for
generic $L$, $M$ into the solution space of a \emph{system}.
Precisely the transition from a higher-order scalar strictly
hyperbolic equation $Lu=0$ in ${\mathbb R}^2$ to an equivalent
first-order characteristic system was used in \cite{Tsarev2005} to
describe a generalization of the Laplace transformation in this
case. For a good definition of general intertwining relations for
linear (probably overdetermined or underdetermined) systems we need
a deeper understanding of the notion of differential transformation
itself since \emph{any} differential mapping of the solution set of
such a general system gives the solution space of another
(overdetermined in general) system unlike the case of scalar
equations $Lu=0$ where description of possible intertwining
relations is not trivial. Probably a generalization of the notion of
$\mathcal{ILT}$ to systems may be of great use. See also \cite{ts09}
for a categorical definition of differential transformations and
factorizations for systems of linear partial differential equations.

So far we did not succeed in representing the important Moutard
transformation \cite{Darboux,gour-l,TT2009} for two-dimensional
stationary Schr\"odinger equation as an $\mathcal{ILT}$. This is a
challenging problem since in the categorical treatment Moutard
transformation is a natural member of the class of
(pseudo)differential transformations in the Serre-Grothendieck
factorcategory of systems (\cite{ts09}). The same similarity of the
Moutard transformation with differential transformations was exposed
in \cite{TT2009} in terms of the skew Ore field of formal fractions
of differential operators.

\section*{Appendix: Intertwining Relations and
Left Least Common Multiples of Linear Partial Differential
Operators}\label{sec-appendix}

It is well known (cf.~for example \cite{Tsarev1996}) that there
always exists the left least common multiple ($lLCM$) for every pair
$L$ and $M$ of linear ordinary differential operators. This is not
always the case for linear partial differential operators $L$ and
$M$. This is related to the algebraic fact that all left (and right)
ideals in $\textbf{F}[D_x]$ are principal ideals, but left ideals in
$\textbf{F}[D_x,D_y]$ are not always principal. In this Appendix we
prove that for many (but not all) examples of intertwining relations
$M_1L=L_1M$ considered in Sect.~\ref{sec-intertw} in fact
$M_1L=L_1M=lLCM(L,M)$.
 Note that in this case if we know the coefficients of the
operators $L$ and $M$ in the intertwining relation (\ref{int1}) we
can find constructively the coefficients of the operators $L_1$ and
$M_1$ from the corresponding system of algebraic equations.
\begin{theorem} Let $L$ and $M$ be linear partial differential
operators in ${\mathbb R}^n$ such that $\mathrm{ord\,} L \geq 1$,
$\mathrm{ord\,} M =\mathrm{ord\,} M_1 =1$ and $L$ is not right
divisible by $M$.
%$rGCD(L,M)=1$.
If $L$ and $M$ satisfy an intertwining relation $M_1L=L_1M$ with
\begin{equation} \label{app-A3}
\mathrm{Sym\,} L = \mathrm{Sym\,} L_1
\end{equation}
then
\begin{equation} \label{app-A4}
lLCM\,(L,M) = M_1L= L_1 M \, .
\end{equation}
\end{theorem}
{\bf Proof.} %\begin{proof}
From the conditions of this theorem we conclude that $L$ and $L_1$
are connected by an intertwining relation (see
Definition~\ref{def-intertw} in Sect.~\ref{sec-intertw}). Hence
$\mathrm{Sym\,} M = \mathrm{Sym\,} M_1$. Suppose that
\begin{equation} \label{app-A5}
 K = P L = Q M
\end{equation}
is some left common multiple of the operators $L$ and $M$. It should
be proved that there exists an operator $G$ such that
\begin{equation} \label{app-A6}
P =G M_1, \qquad  Q = G L_1 \, .
\end{equation}
We prove this by induction on the order of $P$.

\underline{Case 1. \ \ \ \  $ GCD(\mathrm{Sym\,} L ,
 \mathrm{Sym\,} M)=1$}.

Then using (\ref{app-A5}) we see that $\mathrm{Sym\,} P$ should be
divisible by $\mathrm{Sym\,} M$ and  $\mathrm{Sym\,} Q$ should be
divisible by $\mathrm{Sym\,} L$. So we can choose some operator
$G_1$ such that $\mathrm{Sym\,} P = \mathrm{Sym\,} G_1 \cdot
\mathrm{Sym\,} M = \mathrm{Sym\,}  G_1 \cdot \mathrm{Sym\,} M_1$
 and
$\mathrm{Sym\,} Q = \mathrm{Sym\,} G_1 \cdot \mathrm{Sym\,} L =
\mathrm{Sym\,}  G_1 \cdot \mathrm{Sym\,} L_1$. Subtracting from
(\ref{app-A5}) the identity $G_1M_1L = G_1 L_1M$ we obtain
$(P-G_1M_1)L=(Q-G_1L_1)M$ with $\mathrm{ord\,}(P-G_1M_1) <
\mathrm{ord\,} P$,  $\mathrm{ord\,}(Q-G_1L_1) < \mathrm{ord\,} Q$.
So we come to some lower-order left common multiple $K_1=P_1L=Q_1M$.
By induction (if $\mathrm{ord\,} P_1 \geq \mathrm{ord\,} M$,
 $\mathrm{ord\,} Q_1 \geq \mathrm{ord\,} L$) there exists $G_2$ such
that $P_1=G_2M_1$, $Q_1=G_2L_1$, so $P=(G_1+G_2)M_1$,
$Q=(G_1+G_2)L_1$.

If $\mathrm{ord\,} P_1 < \mathrm{ord\,} M$ and $\mathrm{ord\,} Q_1 <
\mathrm{ord\,} L$ then obviously in the case $ GCD(\mathrm{Sym\,} L,
\\ \mathrm{Sym\,} M)=1$ both $P_1$ and $Q_1$ vanish.

\underline{Case 2. \ \ \ \  $ GCD(\mathrm{Sym\,} L ,
 \mathrm{Sym\,} M) \neq 1$}.

Since $\mathrm{ord\,}M =1$ we can choose some operator $S$ such that
$\mathrm{Sym\,} L = \mathrm{Sym\,} M \cdot  \mathrm{Sym\,} S$ and
$\mathrm{Sym\,} L_1 = \mathrm{Sym\,} M_1 \cdot  \mathrm{Sym\,} S$.
Then we have
\begin{equation} \label{app-A7}
L = SM+T, \qquad L_1=M_1S+T_1 \, ,
\end{equation}
with $\mathrm{ord\,}T < \mathrm{ord\,} L$, $\mathrm{ord\,}T_1 <
\mathrm{ord\,} L_1$. From  (\ref{int1}) and (\ref{app-A7}) we obtain
$M_1(SM+T)=(M_1S+T_1)M$ or $M_1T=T_1M$, hence $\mathrm{Sym\,}T =
\mathrm{Sym\,}T_1$.
 \\ If $GCD(\mathrm{Sym\,}M,
\mathrm{Sym\,}T) \neq 1$ we can again simultaneously reduce $T$ and
$T_1$ by $M$ in (\ref{app-A7}) until we obtain  (\ref{app-A7}) with
the condition
\begin{equation} \label{app-A9}
 GCD(\mathrm{Sym\,}M, \mathrm{Sym\,}T) = 1 \, .
\end{equation}
We again take some $LCM(L,M)=K$ and write it in the form
(\ref{app-A5}). Using  (\ref{app-A5}),  (\ref{app-A7}) and
(\ref{app-A9}) we come to $PT=(Q-PS)M$ with $GCD(\mathrm{Sym\,}M,
\mathrm{Sym\,}T) = 1$. As it has been proved in the Case~1 there
exists $G$ such that $P=GM_1$ and $Q-PS=GT_1$. From (\ref{app-A7})
we obtain $T_1=L_1-M_1S$, so $Q-PS=G(L_1-M_1S)$ or $Q=GL_1$.
 \qed
  % \end{proof}

This implies the following proposition:
\begin{proposition} \label{prop-app}
 Let $L$ and $M$ be linear partial differential operators
in ${\mathbb R}^n$ such that $\mathrm{ord\,} L \geq 1$,
$\mathrm{ord\,} M =1$,  $L$ is not right divisible by $M$ and they
satisfy two intertwining relations
\begin{equation} \label{app-A1}
M_1L=L_1M \, ,
\end{equation}
\begin{equation} \label{app-A2}
\widetilde{M}_1L=\widetilde{L}_1M \, ,
\end{equation}
where $\mathrm{Sym\,} L = \mathrm{Sym\,} L_1 = \mathrm{Sym\,}
\widetilde{L}_1$. Then $M_1=\widetilde{M}_1$, $L_1=
\widetilde{L}_1$.
\end{proposition}
{\bf Proof.} %\begin{proof}
Since $lLCM(L,M)$ is unique up to a functional multiplier,
$\widetilde{M}_1 = \phi M_1$, $\widetilde{L}_1 = \phi L_1$ for some
$\phi \in \mathbf{F}$. From the equality $\mathrm{Sym\,} L_1 =
\mathrm{Sym\,} \widetilde{L}_1$ we see that $\phi \equiv 1$.
 \qed
  % \end{proof}

Note that we supposed that $L$ is not divisible by $M$. This is not
always the case---see Sect.~\ref{subsec-gauge} and
Sect.~\ref{subsec-lodo} (Darboux transformations for one-dimensional
Schr\"odinger equations).

%\subsubsection*{Acknowledgments.}
%This paper was written with partial financial support from the KSPU
%grant ``Forming scientific collective \emph{Physics of nano- and
%microstructures}''.


\begin{thebibliography}{99}


\bibitem{Darboux}
Darboux, G.:  Le\c{c}ons sur la th\'eorie g\'en\'erale des surfaces
et les applications g\'eom\'etriques du calcul infinit\'esimal.
t.~2. Gautier-Villard, Paris (1887-1896).


\bibitem{gour-l}
Goursat, E.:  Le\c{c}ons sur l'int\'egration des \'equations aux
d\'eriv\'ees partielles du seconde ordre a deux variables
ind\'ependants. t.~2. Gautier-Villard, Paris, (1898)

\bibitem{Kaptsov} Kaptsov, O.V.: Methods of integration of partial
differential equations (in Russian). Fizmatlit, Moscow (2009)

\bibitem{GLT}
Ganzha, E.I., Loginov, V.M., Tsarev, S.P.:
 Exact solutions of hyperbolic systems of kinetic equations.
Application to Verhulst model with random perturbation.
 Mathematics of Computation.  1, 459--472 (2008)

\bibitem{GanzhaDini}
Ganzha, E.I.: On Laplace and Dini transformations for
multidimensional equations with a decomposable principal symbol.
Programming and Computer Software. 38, 150--155 (2012)

\bibitem{Ore1} Ore, O.: Linear equations in non-commutative fields.
 Ann. Math.  32, 463--477 (1931)

\bibitem{Ore2} Ore, O.:  Theory of non-commutative polynomials.
 Ann. Math. 34, 480--508 (1933)

\bibitem{Tsarev1996}
Tsarev, S.P.:   An algorithm for complete enumeration of all
factorizations of a linear ordinary differential operator. In:
Lakshman, Y.N. (ed.) Proc. ISSAC'1996, pp.~226--231. ACM Press
(1996).


\bibitem{riq}
 Riquier, C.:
Les  syst\`emes d'equations aux deriv\'ees  partielles.
Gautier-Villard, Paris (1910)

\bibitem{jan}
Janet, M.:  Le\c{c}ons sur les syst\`emes d'equations aux deriv\'ees
partielles. Gautier-Villard, Paris. (1929)

\bibitem{schwarz}
Schwarz, F.:
  The Riquier-Janet theory and its applications to nonlinear
evolution equations. Physica D. 11, 243--351 (1984)

\bibitem{She12}
Shemyakova, E.: Laplace transformations as the only degenerate
{D}arboux  transformations of first order. Programming and Computer
Software. 38, 105--108 (2012)


\bibitem{She13b}
Shemyakova, E.: Factorization of Darboux Transformations of
Arbitrary Order for Two-dimensional Schr\"odinger operator,
  \url{http://arxiv.org/abs/1304.7063}  (2013)

\bibitem{TS}
Tsarev S.P., Shemyakova, E.: Differential transformations of
parabolic second-order operators in  the plane. Proc. Steklov Math.
Inst.
%Dedicated  to Academician Sergei Petrovich Novikov on the occasion of his 70th
%  birthday},
266, 219--227 (2009)
%  \url{http://arxiv.org/abs/0811.1492}.

\bibitem{Euler}
Euler, L.: Institutionum calculi integralis. V.~III, Ac. Sc.
Petropoli, St. Petersburg (1770)


\bibitem{LeRoux}
Le Roux, J.:  Extensions de la m\'ethode de Laplace aux \'equations
lin\'eaires aux deriv\'ees   partielles d'ordre sup\'erieur au
second. Bull. Soc. Math. France. 27, 237--262  (1899)
% A digitized copy is obtainable from \url{http://www.numdam.org/}

\bibitem{Petren}
Petr\'en, L.:  Extension de la m\'ethode de Laplace aux \'equations
$\sum_{i=0}^{n-1}A_{1i}\frac{\partial^{i+1}z}{\partial x\partial
y^i}  + \sum_{i=0}^{n}A_{0i}\frac{\partial^{i}z}{\partial y^i} = 0$.
Lund Univ. Arsskrift. 7, 1--166   (1911)

\bibitem{Tsarev2006}
Tsarev, S.P.: On factorization and solution of multidimensional
linear partial differential equations. In: Kotsireas, I., Zima, E.
(eds.) COMPUTER ALGEBRA 2006. Latest Advances in Symbolic
Algorithms, Proc. Waterloo Workshop in Computer Algebra, Canada,
10--12 April 2006, World Scientific. pp.~181-192 (2007)
 \url{http://arxiv.org/abs/cs/0609075}.

\bibitem{Tsarev2005}
Tsarev, S.P.:  Generalized Laplace Transformations and Integration
of Hyperbolic Systems of Linear Partial Differential Equations. In:
Labahn, G. (ed.) Proc. ISSAC'2005.  pp.~325--331. ACM Press (2005)
 \url{http://arxiv.org/abs/cs/0501030}.

\bibitem{Athorne}
Athorne, C.:
  A ${\bf Z}^2 \times {\bf R}^3$ Toda system.
Phys. Lett. A. 206, 162--166  (1995)

\bibitem{ts09}
Tsarev, S.P.:  Factorization of linear differential operators and
systems. In: MacCallum, M.A.H., Mikhailov, A.V. (eds.) Algebraic
Theory of Differential Equations. LMS Lecture Note Series, No.~357.
pp.~111-131 (2009)
 \url{http://arxiv.org/abs/0801.1341}

\bibitem{TT2009}
Taimanov, I.A., Tsarev, S.P.: The Moutard transformation: an
algebraic formalism via pseudodifferential operators and
applications, \url{http://arxiv.org/abs/0906.5141} (2009)


\end{thebibliography}
\end{document}